\documentclass[notitlepage,11pt,reqno]{amsart}
\usepackage[foot]{amsaddr}
\usepackage{amssymb,nicefrac,bm,upgreek,mathtools,verbatim}
\usepackage[final]{hyperref}
\usepackage[mathscr]{eucal}
\usepackage{dsfont}
\usepackage[normalem]{ulem}
\usepackage{amsopn}

\usepackage[margin=1in]{geometry}

\newtheorem{theorem}{Theorem}
\theoremstyle{definition}

\newtheorem{assumption}{Assumption}

\hypersetup{
 colorlinks=true,
 citecolor=mblue,
 linkcolor=mblue,
 urlcolor = mblue,
 anchorcolor = blue,
 frenchlinks=false,
 pdfborder={0 0 0},
 naturalnames=false,
 hypertexnames=false,
 breaklinks}

\usepackage[abbrev,msc-links,nobysame]{amsrefs}

\newcommand{\D}{\mathrm{d}}
\newcommand{\E}{\mathrm{e}}
\newcommand{\Rd}{\mathbb{R}^d}
\newcommand{\RR}{\mathbb{R}}
\DeclareMathOperator{\Exp}{\mathrm{E}}

\newcommand{\Ind}{\bm{1}} 

\newcommand{\Act}{\mathbb{U}}

\newcommand{\Cc}{{C}} 
\newcommand{\cG}{{\mathcal{G}}} 
\newcommand{\sG}{{\mathscr{G}}} 
\newcommand{\cH}{{\mathcal{H}}} 
\newcommand{\eom}{{\mathcal{M}}} 
\newcommand{\varphis}{{\varphi^{}_{\mspace{-2mu}*}}}
\newcommand{\Phis}{\Phi_{\mspace{-2mu}*}}
\newcommand{\cPs}{{\cP^{}_{\mspace{-3mu}*}}}
\newcommand{\cPo}{{\cP^{}_{\mspace{-3mu}\circ}}}

\newcommand{\cZ}{{\mathcal{Z}}}
\newcommand{\cP}{{\mathcal{P}}} 
\newcommand{\sA}{\widetilde{\mathcal{A}}} 
\newcommand{\cA}{\mathcal{A}}
\newcommand{\cL}{\mathcal{L}}

\newcommand{\cK}{\mathcal{K}}

\newcommand{\Lyap}{\mathscr{V}} 

\newcommand{\df}{\coloneqq}
\newcommand{\transp}{^{\mathsf{T}}}

\newcommand{\grad}{\nabla}
\DeclareMathOperator*{\trace}{trace}
\newcommand{\sorder}{{\mathfrak{o}}} 

\newcommand{\Sobl}{{\mathscr W}_{\text{loc}}} 
\newcommand{\Lpl}{L_{\text{loc}}} 

\newcommand{\abs}[1]{\lvert#1\rvert}
\newcommand{\norm}[1]{\lVert#1\rVert}
\newcommand{\babs}[1]{\bigl\lvert#1\bigr\rvert}

\usepackage{color}
\definecolor{dmagenta}{rgb}{.4,.1,.5}
\definecolor{dblue}{rgb}{.0,.0,.5}
\definecolor{mblue}{rgb}{.0,.0,.55}
\definecolor{ddblue}{rgb}{.0,.0,.4}
\definecolor{dred}{rgb}{.7,.0,.0}
\definecolor{dgreen}{rgb}{.0,.5,.0}
\definecolor{Eeom}{rgb}{.0,.0,.5}

\newcommand{\ttl}{\Large `Controlled' versions of the Collatz--Wielandt\\[3pt]
 and Donsker--Varadhan formulae}
\begin{document}
\title[`Controlled' versions of the Collatz--Wielandt
 and Donsker--Varadhan formulae]{\ttl}

\author[Ari Arapostathis]{Ari Arapostathis$^*$}
\address{$^*$Department of Electrical and Computer Engineering,
The University of Texas at Austin,
EER~7.824, Austin, TX~~78712}
\email{ari@utexas.edu}

\author[Vivek S. Borkar]{Vivek S. Borkar$^\ddag$}
\address{$^\ddag$Department of Electrical Engineering,
Indian Institute of Technology, Powai, Mumbai 400076, India}
\email{borkar.vs@gmail.com}



\begin{abstract}
This is an overview of the work of the authors and their collaborators on
the characterization
of risk sensitive costs and rewards in terms of an abstract Collatz--Wielandt formula
and in case of rewards, also a controlled version of the Donsker--Varadhan formula.
For the finite state and action case, this leads to useful linear and dynamic
programming formulations in the reducible case.

\end{abstract}

\maketitle

\section{Introduction}

This short article is an overview of the work of authors and their collaborators on
a somewhat novel perspective of the risk-sensitive control problem on infinite time
horizon that aims to optimize the asymptotic growth rate of a mean exponentiated total
reward, resp., cost.
The viewpoint taken here is based on the fact that the dynamic programming principle
for this problem essentially reduces it to an eigenvalue problem
seeking the principal eigenvalue and eigenvector for a monotone positively $1$-homogeneous
operator. This allows us to exploit the existing generalized Perron--Frobenius
(or Krein--Rutman) theory which leads to some explicit expressions for the optimal
growth rate. The first is the abstract
Collatz--Wielandt formula which can be shown to hold for both cost minimization and
reward maximization problems, though we have not exhausted all the cases in our work.
The second is a variational formula for the
principal eigenvalue that generalizes the Donsker--Varadhan formula for the same in
the linear case. This seems workable only for the reward maximization problem.

We first consider the discrete time case based on the results of \cite{Ananth}
in the next two sections, followed by those for reflected diffusions in a bounded domain,
based on \cite{ABK}, in section~\ref{S4}. We then sketch, in section~\ref{S5}, the very
recent and highly nontrivial extensions to diffusions on the whole space developed in
\cite{AriAnup} and \cite{AABK}.
Finally, we recall in section~\ref{S6} some developments in the simple finite
state-action set up from \cite{CDC}, where the aforementioned development allows
us to derive the dynamic programming equations for risk-sensitive reward process
in the reducible case. Section~7 concludes by highlighting some future directions.

\section{Discrete time problems}

The celebrated Courant--Fisher formula for the principal eigenvalue of a positive
definite symmetric matrix $A \in \RR^{d\times d}$ is
\begin{equation*}
\lambda = \max_{0 \neq x \in \Rd}\frac{x\transp Ax}{x\transp x}.
\end{equation*}
Consider an irreducible nonnegative matrix $Q \in \RR^{d\times d}$.
The Perron--Frobenius theorem guarantees a positive principal eigenvalue with an
associated positive eigenvector for $Q$.
Is there a counterpart of the Courant--Fisher formula for this eigenvalue?

The answer is a resounding `YES'! It is the Collatz-Wielandt formula for the
principal eigenvalue of an irreducible nonnegative matrix
$Q = \boldsymbol[q(i,j)\boldsymbol] \in \RR^{d\times d}$, stated as
(see \cite{Meyer} Chapter~8):
\begin{align*}
 \lambda &\,=\, \sup_{x = [x_1, \cdots, x_d]\transp,\, x_i \geq 0 \ \forall i} \;
 \min_{i \colon x_i > 0} \; \left(\frac{(Qx)_i}{x_i}\right) \\[5pt]
 &\,=\, \inf_{x = [x_1, \cdots, x_d]\transp,\, x_i > 0 \ \forall i} \;
 \max_{i \colon x_i > 0}
 \; \left(\frac{(Qx)_i}{x_i}\right).
 \end{align*}
An alternative characterization can be given as follows. Write
\begin{align*}
Q &\,=\, \varGamma P\,,\\
\intertext{where}
\varGamma &\,\df\, \ \text{diag}(\kappa_1, \dotsc, \kappa_d)\,,\ \ \kappa_i>0\,,\\
p(i,j) &\,\df\, \nicefrac{q(i,j)}{\kappa_i}\,, \quad 1 \leq i,j \leq d\,, \\
P &\,\df\, \boldsymbol[p(j\,|\,i)\boldsymbol]\,,
\end{align*}
with $P$ a stochastic matrix.
In other words, we have pulled out the row sums $\{\kappa_i\}$ of $Q$ into a diagonal
matrix $\varGamma$ so that what is left is a stochastic matrix $P$.
Also define
\begin{equation*}
\sG_0 \,\df\, \bigl\{(\pi, \tilde{P}) \colon \pi \text{\ is a stationary probability
for the stochastic matrix\ }
\tilde{P} = \boldsymbol[\tilde{p}(j|i)\boldsymbol]\bigr\}\,.
\end{equation*}
Then the following representation holds \cite{Dembo}:
\begin{equation*}
\log\lambda \,=\, \sup_{(\pi, \tilde{P}) \,\in\, \sG_0}
\left(\sum_i\pi(i)\bigl[\kappa_i - D\bigl(\tilde{p}(\cdot\,|\, i ) \|\,
p( \cdot\,|\, i )\bigr)\bigr]\right),
\end{equation*}
where $D( \cdot \,\|\, \cdot)$ denotes the Kullback--Leibler divergence or relative entropy.
This is the finite state counterpart of the Donsker--Varadhan formula \cite{DoVa}
for the principal eigenvalue of a nonnegative matrix.

As is well known, the infinite dimensional generalization of the Perron--Frobenius
theorem is given by the Krein--Rutman theorem \cite{Krein,Pagter}.
There are also nonlinear variants of it.
Let
\begin{enumerate}
\item $B$ be a Banach space with a `positive cone' $K$ such that $K - K$ is dense in $B$,

\item $T\colon B \mapsto B$ be a compact order preserving
(i.e., $f \geq g \Longrightarrow Tf \geq Tg$), strictly increasing
(i.e., $f >g \Longrightarrow Tf > Tg)$, strongly positive
(i.e., maps nonzero elements of $K$ to its interior),
positively $1$-homogeneous (i.e., $T(af) = aTf$ for all $a > 0$) operator.
\end{enumerate}
A nonlinear variant of the Krein--Rutman theorem \cite{Ogiwara} then asserts that under
some technical hypotheses, a unique positive principal eigenvalue and a corresponding
unique (up to a scalar multiple) positive eigenvector for $T$ exist.

Our interest is in the following \emph{nonlinear} scenario arising in
risk-sensitive control: Consider
\begin{itemize}
\item a controlled Markov chain $\{X_n\}$ on a compact metric state space $S$;

\item an associated control process $\{Z_n\}$ in a compact metric control space $U$;

\item a \emph{per stage reward function} $r \colon S\times U\times S \mapsto \RR$
such that $r \in C(S\times U\times S)$;

\item a controlled transition kernel $p(\D y\,|\,x,u)$ with \textbf{full support},
such that for all Borel $A \subset S$,
\begin{equation}\label{CMarkov}
\begin{aligned}
P(X_{n+1} \in A \,|\, X_m, Z_m, m \leq n) &\,=\, P(X_{n+1} \in A \,|\, X_n, Z_n) \\
&\,=\, p(A\,|\, X_n,Z_n)\,.
\end{aligned}
\end{equation}
This is called the \emph{controlled Markov property} and the controls for which this
holds are said to be \emph{admissible}. The maps
\begin{equation*}
(x,u) \mapsto \int f(y)p(\D y\,|\,x,u), \quad f \in C(S), \ \|f\| \leq 1,
\end{equation*}
are assumed to be equicontinuous.
\end{itemize}

The control problem is to maximize the asymptotic growth rate of the exponential reward:
\begin{equation*}
\lambda \,\df\, \sup_{x \in S}\,\sup_{\{Z_m\}}\,\liminf_{N\uparrow\infty}\frac{1}{N}
\log \Exp\left[\E^{\sum_{m=0}^{N-1}r(X_m, Z_m, X_{m+1})}\Bigm|X_0 = x\right].
\end{equation*}
The second supremum in this definition is over all admissible controls.
We allow \emph{relaxed} (i.e., probability measure valued) controls $\{\mu_n\}$
taking values in $\cP(S)$, in which case \eqref{CMarkov} gets replaced by
\begin{equation*}
\begin{aligned}
P(X_{n+1} \in A \,|\, X_m, \mu_m, \ m \leq n) &\,=\, P(X_{n+1} \in A\,|\, X_n, \mu_n) \\
&\,=\, \int p(A \,|\, X_n, z)\mu_n(dz), \ n \geq 0\,.
\end{aligned}
\end{equation*}
Define
$$Tf(x) \,\df\, \sup_{\phi : S \mapsto \cP(U) \text{\ measurable}}\;
\iint p(\D y\,|\,x,u)\phi(\D u\,|\,x)\E^{r(x,u,y)}f(y)\,.$$
This is a compact, order preserving, strictly increasing, strongly positive,
positively $1$-homogeneous operator.

Using the nonlinear variant of the Krein--Rutman theorem stated above,
this leads to an abstract Collatz-Wielandt formula \cite{Ananth}:

\begin{theorem}
There exist $\rho > 0, \psi \in \mathrm{int}(C^+(S))$ such that $T\psi = \rho\psi$ and
\begin{equation*}
\begin{aligned}
\rho &\,=\, \inf_{f \,\in \,\mathrm{int}(C^+(S))} \; \sup_{\eom^+(S)} \;
\frac{\int Tfd\mu}{\int fd\mu}\\
&\,=\, \sup_{f \,\in \,\mathrm{int}(C^+(S))} \; \inf_{\eom^+(S)} \;
\frac{\int Tfd\mu}{\int fd\mu}\,.
\end{aligned}
\end{equation*}
Also, $\log\rho$ is the optimal reward for the risk-sensitive control problem.
\end{theorem}

\section{Variational Formula}

We now state a variational formula for the principal eigenvalue \cite{Ananth}.
Let $\sG$ denote
the set of probability measures
\begin{equation*}
\eta(\D x,\D u,\D y) \in \cP(S\times U\times S)
\end{equation*}
which disintegrate as
\begin{equation*}
\eta(\D x,\D u,\D y) \,=\, \eta_0(\D x)\eta_1(\D u\,|\,x)\eta_2(\D y\,|\,x,u)\,,
\end{equation*}
such that $\eta_0$ is invariant under the transition kernel
\begin{equation*}
\int_U \eta_2(\D y\,|\,x,u)\eta_1(\D u\,|\,x)\,.
\end{equation*}
These are the so called `ergodic occupation measures' for discrete time control problems.

\begin{theorem} Under the above hypotheses,
\begin{equation*}
\log\rho \,=\, \sup_{\eta \in \sG} \biggl( \iint\eta_0(\D x)\eta_1(\D u\,|\,x)
\biggl[\int r(x,u,y)\eta_2(\D y\,|\,x,u) 
- \ D\bigl(\eta_2(\D y\,|\,x,u)\|\, p(\D y\,|\,x,u)\bigr)\biggr] \biggr).
\end{equation*}
\end{theorem}

\bigskip

This can be viewed as a controlled version of the Donsker--Varadhan formula.
The hypotheses above can be relaxed to:
\begin{enumerate}
\item Range$(r) = [-\infty, \infty)$ with $\E^r \in C(S\times U\times S)$;

\item $p(\D y\,|\,x,u)$ need not have full support.
\end{enumerate}
The formula then is the same as before, the difference is that under the previous,
stronger set of conditions, the supremum over $x \in S$ in the definition of $\lambda$
was redundant, it is no longer so.
The extension proceeds via an approximation argument that approximates the given
transition kernel by a sequence of transition kernels for which our original hypotheses hold.

We thus have an equivalent concave maximization problem, in fact a linear program,
as opposed to a `team' problem one would obtain from the usual `log transformation' as in,
e.g., \cite{Flem}.
Furthermore, if $\rho(\varphi)$ denotes the asymptotic growth rate for a randomized Markov
control $\varphi$, then it can be shown that $\rho = \max_{\varphi} \rho(\varphi)$,
implying the sufficiency of randomized Markov controls.

Some applications worth noting are \cite{Ananth}:
\begin{enumerate}

\item Growth rate of the number of directed paths in a graph.
This requires $-\infty$ as a possible reward to account for the absence of edges.

\item Portfolio optimization in the framework of \cite{Bielecki}.

\item Problem of minimizing the exit rate from a domain.
\end{enumerate}

\section{Reflected diffusions}\label{S4}

Analogous results hold for reflected diffusions in a compact domain with smooth boundary.
These are described by the stochastic differential equation
\begin{equation}\label{E-sderefl}
\begin{aligned}
 \D X(t) &\,=\, b\bigl(X_t, U_t\bigr)\,\D t + \upsigma\bigl(X_t\bigr)\,\D W_t
 - \gamma(X_t)\,\D\xi_t\,, \\[2pt]
 \D \xi(t) &\,=\, \bm1\{X_t \in \partial Q\}\,\D\xi_t\,,
\end{aligned}
\end{equation}
 for $t \geq 0$. Here:
 \begin{enumerate}

 \item $Q$ is an open connected and bounded set with $C^3$ boundary $\partial Q$;

\item $\{W_t\}_{t\ge0}$ is a standard $d$-dimensional Wiener process;
\item
the control $\{U_t\}_{t\ge0}$ lives in a metrizable compact action space $\Act$ and is
non-anticipative, i.e., for $t > s$, $W(t) - W(s)$ is independent of
$X_0; W_y, U_y, y \leq s$;

 \item $b$ is continuous, and $x\mapsto b(x, u)$ is Lipschitz uniformly in $u$;

 \item $\upsigma$ is $C^{1,\beta_0}$ and uniformly non-degenerate;

 \item $\gamma_i(x) = \upsigma(x)\upsigma(x)\transp\eta(x)$ where $\eta(x)$
 is the unit outward normal
 on $\partial Q$.

 \end{enumerate}

 In contrast to the preceding section, we first consider the cost minimization problem
 to highlight the differences with the reward maximization problem.
 Unlike the classical cost/reward criteria such as discounted and average cost/reward,
 the risk-sensitive cost and reward problems are not rendered equivalent by a mere
 sign flip, and the differences are stark.
 For cost minimization, the control problem is to minimize
 \begin{equation*}
 \lim_{t\uparrow\infty}\,
 \frac{1}{t}\,\log \Exp\left[\E^{\int_0^tr(X_s, U_s)\,\D s}\right],
 \end{equation*}
 where $r$ is continuous.

 The corresponding `Nisio semigroup' is defined as follows. For $t \geq 0$, let
\begin{equation*}
 S_tf(x) \,\df\, \inf_{\{U_t\}_{t\ge0}}\,\Exp_x\left[\E^{\int_0^t
 r(X_s, U_s)\,\D s}f(X_t)\right].
 \end{equation*}
 Then $S_t \colon C(\bar{Q}) \mapsto C(\bar{Q})$ is a semigroup of strongly continuous,
 bounded Lipschitz, monotone, superadditive, positively 1-homogeneous, strongly positive
 operators with infinitesimal generator $\cG$ defined by
 \begin{equation}\label{E-cG}
\cG f(x) \,\df\, \frac{1}{2}\mbox{tr}\left(\upsigma(x)\upsigma\transp(x)\nabla^2f(x)\right) +
\min_{u\in\Act}\,\Bigl[\langle b(x,u)\,, \nabla f(x)\rangle + r(x,u)f(x)\Bigr]\,.
 \end{equation}

 Let
 \begin{equation*}
 C^2_{\gamma,+}(\bar{Q}) \,\df\,
 \bigl\{f\colon\bar{Q} \mapsto [0, \infty) \colon f\in C^2(\Bar{Q}),
 \ \langle\nabla f(x), \gamma(x)\rangle = 0 \text{\ for\ }x \in \partial Q\bigr\}\,.
 \end{equation*}
 As in the discrete case, the nonlinear Krein--Rutman theorem then leads to:
 There exists a unique pair $(\rho, \varphi) \in \RR\times C^2_{\gamma, +}(\bar{Q})$
 satisfying $\|\varphi\|_{0,\bar{Q}} = 1$ such that
 \begin{equation*}
 S_t\varphi = \E^{\rho t}\varphi\,.
 \end{equation*}
 This solves
 \begin{equation*}
\cG \varphi(x) \,=\, \rho\varphi(x)\,,\ \ x \in Q, \qquad\text{and\ \ }
\langle \nabla\varphi(x), \gamma(x)\rangle \,=\, 0\,, \ \ x \in \partial Q\,.
 \end{equation*}

 The abstract Collatz-Wielandt formula for this problem is
\begin{equation*}
\begin{aligned}
 \rho &\,=\, \inf_{f \in C^2_{\gamma,+}(\bar{Q}), f > 0}\;\sup_{\nu \in \cP(\bar{Q})}\;
 \int_{\bar{Q}}\frac{\cG f}{f}\D\nu \\
 &\,=\, \sup_{f \in C^2_{\gamma,+}(\bar{Q}), f > 0}\;\inf_{\nu \in \cP(\bar{Q})}\;
 \int_{\bar{Q}}\frac{\cG f}{f}\D\nu\,.
\end{aligned}
\end{equation*}
 In the uncontrolled case, the first formula above is the convex dual of the
 Donsker--Varadhan formula for the principal eigenvalue of $\cG$:
 \begin{equation*}
 \rho \,=\, \sup_{\nu \in \cP(\bar{Q})}\left(\int_{\bar{Q}}r(x)\nu(\D x) - I(\nu)\right)\,,
 \end{equation*}
 where
 \begin{equation*}
 I(\nu) \,\df\,
 \inf_{f \in C^2_{\gamma,+}(\bar{Q}), f > 0}\;\int_{\bar{Q}}\left(\frac{\cG f}{f}\right)
 \D\nu\,.
 \end{equation*}

For the risk-sensitive reward problem, the same abstract Collatz-Wielandt formula holds,
except that the definition of the operator $\cG$ now has a `$\max$' in place of the `$\min$'.
But as in the discrete time case, one can go a step further and have a
variational formulation. Let
\begin{equation*}
R(x,u,w) \,\df\, r(x,u) - \frac{1}{2}\abs{\upsigma\transp(x)w}^2\,,
\quad (x,u,w)\in\Bar{Q}\times\Act\times\Rd\,,
\end{equation*}
and
\begin{equation*}
 \eom \,\df\, \biggl\{\mu \in \cP(\bar{Q}\times U\times \Rd) \,\colon 
 \int_{\bar{Q}\times U\times \Rd}\cA f(x,u,w)
 \mu(\D x,\D u,\D w) = 0\ \ \forall\, f \in C^2(Q)\cap C_\gamma(\Bar{Q})\biggr\}\,,
\end{equation*}
with
\begin{equation}
\cA f (x,u,w)\,\df\, \frac{1}{2}\mbox{tr}
\left(\upsigma(x)\upsigma\transp(x)\nabla^2f(x)\right) +
 \bigl\langle b(x,u) + \upsigma(x)\upsigma\transp(x)w, \nabla f(x)\bigr\rangle
 \label{E-cA}
\end{equation}
for $f \in C^2(Q)\cap C(\bar{Q})$.
Recall the definition of an `\emph{ergodic occupation measure}' \cite{ABG}.
For a stochastic differential equation as in \eqref{E-sderefl}, but with the drift $b$
replaced with
$b(x,u) + \upsigma(x)\upsigma\transp(x)w$, and $w$ taking values in some
compact metrizable space,
it is the time-$t$ marginal of a stationary state-control process
$\bigl(X_t, v(X_t), w(X_t)\bigr)$,
perforce independent of $t$.
Thus, in the case the parameter $w$ lives in a compact space,
by a standard characterization of ergodic occupation measures (\emph{ibid.}),
$\eom$ is precisely the set thereof
for controlled diffusions whose (controlled) extended generator is $\cA$.
This however is not necessarily the case if $w$ lives in $\Rd$.
An example to keep in mind is the one-dimensional
stochastic differential equation
\begin{equation*}
\D{X}_{t} = \bigl(\E^{\nicefrac{X_{t}^{2}}{2}}- X_{t}\bigr)\,\D{t} +\sqrt{2}\, \D{W}_{t}\,.
\end{equation*}
It is straightforward to verify that the standard Gaussian density satisfies the
Fokker--Planck equation. However, the diffusion is not even regular, so
it does not have an invariant probability measure.
Therefore, we refer to $\eom$ as the set of \emph{infinitesimal ergodic occupation
measures}.
The variational formula for this model is
$$\rho \,=\, \sup_{\mu \in \eom}\int_{\bar{Q}\times U\times \Rd}R(x,u,w)
\mu(\D x,\D u,\D w)\,.$$
This result is from \cite{AABK}.

An analogous abstract Collatz--Wielandt formula for the risk-sensitive
\emph{cost minimization} problem was derived in \cite{ABK}.
We have not derived a corresponding variational formula.
Even if one were to do so, it is clear that it will be a `sup-inf / inf-sup' formula
rather than a pure maximization problem.
This is already known through a different route: it forms the basis of the approach
initiated by \cite{Flem} and followed by many, in which the the Hamilton--Jacobi--Bellman
equation for the risk-sensitive cost minimization problem is converted to an Isaacs
equation for an ergodic payoff zero sum stochastic differential game.
The aforementioned expression then is simply the value of this game.
Going by pure analogy, for the reward maximization problem, one would expect this route
to yield a stochastic \emph{team} problem wherein the two agents seek to maximize a
common payoff, but \emph{non-cooperatively}, i.e., without either of them having
knowledge of the other person's decision.
What this translates into is that under the corresponding ergodic occupation measure,
the two control actions are conditionally independent given the state.
The set of such measures is non-convex.
What we have achieved instead is a single concave programming problem, which is a
significant simplification from the point of view of developing computational schemes
for the problem.
This also brings to the fore the difference between reward maximization and cost
minimization in risk-sensitive control.

\section{Diffusions on the whole space}\label{S5}

Here we consider a controlled diffusion in
$\Rd$ of the form
\begin{equation*}
\D X_t \,=\, b(X_t,U_t)\,\D t + \upsigma (X_t)\,\D W_t\,,
\end{equation*}
where
\begin{enumerate}
\item $W$ is a standard $d$-dimensional Brownian motion;
\item
the control $U_t$ lives in a metrizable compact action space $\Act$ and is
non-anticipative, i.e., for $t > s$, $W(t) - W(s)$ is
independent of $X_0; W_y, U_y, y \leq s$;
\item
 $b(x,u)$ is continuous and locally Lipschitz continuous in $x$ uniformly in $u\in\Act$;
\item
$\upsigma$ is locally Lipschitz continuous and locally nondegenerate;
\item
$b$ and $\upsigma$ have at most affine growth in $x$.
\end{enumerate}
Without loss of generality, we may take $U_t$ to be adapted to the increasing
$\sigma$-fields generated by $\{X_t, t \geq 0\}$.
Then these hypotheses guarantee the existence of a unique weak solution for
any admissible control $\{U_t\}_{t\ge0}$ (\cite{ABG}, Chapter~2).

As before, we let $r(x,u)$ be a continuous running reward function, which is
locally Lipschitz in $x$ uniformly in $u$, and is also bounded from above in
$\Rd$.
We define the
\emph{optimal risk-sensitive value} $J^*$ by
\begin{equation*}
J^* \,\df\, \sup_{\{U_t\}_{t\ge0}}\;\liminf_{T\to\infty}\, \frac{1}{T}\,
\log \Exp \Bigl[\E^{\int^T_0 r(X_t,U_t)\,\D t} \Bigr]\,,
\end{equation*}
where the supremum is over all admissible controls.

Consider the extremal operator
\begin{equation*}
\widehat\cG f(x) \,\df\, \frac{1}{2}\trace\left(a(x)\nabla^{2}f(x)\right)
+ \max_{u\in\Act}\, \Bigl[\bigl\langle b(x,u),
\grad f(x)\bigr\rangle + r(x,u) f(x)\Bigr]
\end{equation*}
for $f\in\Cc^2(\Rd)$.
The \emph{generalized principal eigenvalue} of $\widehat\cG$ is defined by
\begin{equation}\label{E-lamstr}
\lambda_*(\widehat\cG)\,\df\,\inf\,\Bigl\{\lambda\in\RR\,
\colon \exists\, \phi\in\Sobl^{2,d}(\Rd),\ \varphi>0, \
\widehat\cG\phi -\lambda\phi\le 0 \text{\ a.e.\ in\ } \Rd\Bigr\}\,,
\end{equation}
where $\Sobl^{2,d}(\Rd)$ denotes the local Sobolev space of
functions on $\Rd$ whose generalized
derivatives up to order $2$ are in $\Lpl^{d}(\Rd)$, equipped with its natural
semi-norms.
We assume that $r-\lambda_*$ is negative and
bounded from above away from zero on the complement of some compact set.
This is always satisfied if $-r$ is an inf-compact function, that is
the sublevel sets $\{-r \le c\}$ are compact (or empty) in $\Rd\times\Act$
for each $c\in\RR$, or if $r$ is a positive function vanishing at infinity
and the process $\{X_t\}_{t\ge0}$ is recurrent under some
stationary Markov control.
Then there exists a unique positive $\Phis\in\Cc^2(\Rd)$ normalized as
$\Phis(0)=1$ which solves $\widehat\cG\Phis = \lambda_*\Phis$.
In other words, the eigenvalue $\lambda_*=\lambda_*(\widehat\cG)$ is simple.
Let $\varphis\df\log\Phis$.
As shown in \cite{AABK}, the function
\begin{equation*}
\cH(x)\,\df\,\frac{1}{2}\,\babs{\upsigma\transp(x)\grad \varphis(x)}^2\,,
\quad x\in\Rd
\end{equation*}
is an \emph{infinitesimal relative entropy rate}.

We let $\cZ \df \Rd\times\Act\times\Rd$, and
use the single variable $z=(x,u,w)\in\cZ$.
Let $\cP(\cZ)$ denote the set of probability measures
on the Borel $\sigma$-algebra of $\cZ$, and
$\eom_A$ denote the set of infinitesimal ergodic occupation measures
for the operator $\cA$ in \eqref{E-cA} defined for $f\in\Cc^2(\Rd)$,
which here can be written as
\begin{equation*}
\eom_{\cA}\,\df\,\biggl\{ \mu\in \cP(\cZ)\,\colon
\int_{\cZ} \cA f(z)\,\mu(\D{z}) \,=\, 0\quad \forall\, f\in\Cc^{2}_c(\Rd)\biggr\}\,,
\end{equation*}
where $\Cc^2_c(\Rd)$ is the class of functions
in $\Cc^2(\Rd)$ which have compact support.
Recall the definition $R(x,u,w)\df r(x,u) - \frac{1}{2}\abs{\upsigma\transp(x)w}^2$
in Section~\ref{S4}.
We also define
\begin{equation*}
\begin{aligned}
\cPs(\cZ) &\,\df\,
\biggl\{\mu\in \cP(\cZ)\,\colon
\int_{\cZ} \cH(x)\,\mu(\D{x},\D{u},\D{w}) <\infty\biggr\}\,,\\
\cPo(\cZ) &\,\df\,
\biggl\{ \mu\in \cP(\cZ)\,\colon
\int_{\cZ} R(z)\,\mu(\D{z}) > -\infty\biggr\}\,.
\end{aligned}
\end{equation*}

The following is a summary of the main results in \cite[Section~4]{AABK}.

\begin{theorem}
We have
\begin{equation*}
\begin{aligned}
J^* \,=\, \lambda_*(\widehat\cG) &\,=\, \sup_{\mu\in\cPs(\cZ)}\,
\inf_{g \in \Cc^2_c(\Rd)}\,\int_{\cZ} \bigl(\cA g(z)+R(z)\bigr)\,\mu(\D{z})\\
&
\,=\, \max_{\mu\in\eom_{\sA}\cap\cPs(\cZ)}\,\int_{\cZ} R(z)\,\mu(\D{z})\,.
\end{aligned}
\end{equation*}
Suppose that the diffusion matrix
$a$ is bounded and uniformly elliptic, and either
$-r$ is inf-compact, or $\langle b,x\rangle^-$ has subquadratic growth,
or $\frac{\abs{b}^2}{1+\abs{r}}$ is bounded.
Then $\eom_{\cA}\cap\cPo(\cZ)\subset\cPs(\cZ)$, and
$\cPs(\cZ)$ may be replaced by $\cP(\cZ)$ in the variational formula above.
If, in addition, $\frac{\cH}{1+\abs{\varphis}}$ is bounded, then
\begin{equation*}
J^* \,=\, \lambda_*(\widehat\cG) \,=\, \inf_{g \in \Cc^2_c(\Rd)}\, \sup_{\mu\in\cP(\cZ)}\,
\int_{\cZ} \bigl(\cA g(z)+R(z)\bigr)\,\mu(\D{z}) \,.
\end{equation*}
\end{theorem}

We continue with the Collatz--Wielandt formula in $\Rd$ for the risk-sensitive
cost minimization problem. This is studied in \cite{AriAnup}.
Here, we have
a running cost $r(x,u)$ which is bounded from below in $\Rd\times\Act$,
and is locally Lipschitz in $x$ uniformly in $u$.
The assumptions on $b$ and $\upsigma$ are as stated in the beginning of
the section, except that we may replace the affine growth assumption with
the more general condition
\begin{equation*}
\sup_{u\in\Act}\; \langle b(x,u),x\rangle^{+} + \norm{\upsigma(x)}^{2}\,\le\,C_0
\bigl(1 + \abs{x}^{2}\bigr) \qquad \forall\, x\in\RR^{d}\,,
\end{equation*}
for some constant $C_0>0$.
The risk-sensitive optimal value $\Lambda^*$ is defined by
\begin{equation*}
\Lambda^*\,\df\, \inf_{\{U_t\}_{t\ge0}}\; \limsup_{T\to\infty}\,\frac{1}{T}\,
\log\Exp\Bigl[\E^{\int_0^T r(X_s, U_s)\, \D{s}}\Bigr]\,.
\end{equation*}
The operator $\cG$ here is as in \eqref{E-cG} but for $f\in\Cc^2(\Rd)$,
and we let the generalized principal eigenvalue $\lambda_*(\cG)$ be defined
as in \eqref{E-lamstr}.

The running cost does not have any structural properties
that penalize unstable behavior
such as near-monotonicity or inf-compactness, so uniform ergodicity
for the controlled process needs to be assumed.
Let
\begin{equation*}
\cL f(x,u) \,\df\,
\frac{1}{2}\mbox{tr}\left(\upsigma(x)\upsigma\transp(x)\nabla^2f(x)\right)
+ \bigl\langle b(x,u), \nabla f(x)\bigr\rangle\,.
\end{equation*}
We consider the following hypothesis.

\begin{assumption}\label{A2.1}
The following hold.
\begin{itemize}
\item[(i)] There exists an inf-compact function $\ell\in\Cc(\Rd)$, and a positive
function $\Lyap\in\Sobl^{2, d}(\Rd)$, satisfying $\inf_{\Rd}\Lyap>0$, such that
\begin{equation}\label{EA2.1A}
\sup_{u\in\Act}\,\cL\Lyap(x,u) \,\le\, \kappa_1 \Ind_{\cK}(x)
-\ell(x) \Lyap(x) \quad \forall\,x\in\Rd\,,
\end{equation}
for some constant $\kappa_1$ and a compact set $\cK$.
\item[(ii)]
The function $x\mapsto \beta\ell(x)-\max_{u\in\Act}\, r(x,u)$ is inf-compact
for some $\beta\in(0, 1)$.
\end{itemize}
\end{assumption}

As noted in \cite{ABS-19}, the Foster--Lyapunov equation
in \eqref{EA2.1A} cannot in general be satisfied
 for diffusions with bounded $a$
and $b$. Therefore, to treat this case,
we consider an alternate set of conditions.

\begin{assumption}\label{A2.2}
The following hold.
\begin{itemize}
\item[(i)]
There exists a positive function $\Lyap\in\Sobl^{2,d}(\Rd)$,
satisfying $\inf_{\Rd} \Lyap > 0$, constants $\kappa_1$ and $\gamma>0$,
and a compact set $\cK$ such that
\begin{equation*}
\sup_{u\in\Act}\,\cL\Lyap(x,u) \le \kappa_1 \Ind_{\cK}(x)
-\gamma \Lyap(x) \quad \forall\,x\in\Rd\,.
\end{equation*}
\item[(ii)]
$\norm{r^-}_\infty+\limsup_{\abs{x}\to\infty}\, \max_{u\in\Act}\,r(x,u)
<\gamma$.
\end{itemize}
\end{assumption}

Let $\sorder(\Lyap)$ denote the class of continuous functions $f$ that grow slower
than $\Lyap$, that is, $\frac{\abs{f(x)}}{\Lyap(x)}\to0$ as $\abs{x}\to\infty$.
We quote the following result from \cite{ABS-19}.

\begin{theorem}
Grant either Assumption~\ref{A2.1}, or~\ref{A2.2}.
Then
\begin{equation}\label{ET2.4A}
\begin{aligned}
\Lambda^*\,=\,\lambda_*(\cG) &\,=\, \sup_{f\in \Cc^{2,+}(\Rd)\cap\sorder(\Lyap)}\;
\inf_{\mu\in\cP(\Rd)}\; \int_{\Rd}\frac{\cG f}{f}\, \D{\mu}\\
&\,=\, \inf_{f\in \Cc^{2,+}(\Rd)}\;
\sup_{\mu\in\cP(\Rd)}\; \int_{\Rd}\frac{\cG f}{f}\, \D{\mu}\,,
\end{aligned}
\end{equation}
where $\Cc^{2,+}(\Rd)$ denotes the set of positive functions in $\Cc^2(\Rd)$.
\end{theorem}

We should remark here that the class of test functions $f$ in the first
representation formula in \eqref{ET2.4A} cannot,
in general, be enlarged to $\Cc^{2,+}(\Rd)$.

It is also interesting to consider the substitution $f = e^{\psi}$. Then \eqref{ET2.4A}
transforms to
\begin{align*}
\lambda_*(\cG) &\,=\,
\sup_{\psi\in \Cc^{2,+}(\Rd)\cap\sorder(\log\Lyap)}\; \inf_{\mu\in\cP(\Rd)}\;
\int_{\Rd}F[\psi](x)\,\mu(\D{x})\\
&\,=\, \inf_{\psi\in \Cc^{2,+}(\Rd)}\; \sup_{\mu\in\cP(\Rd)}\;
\int_{\Rd}F[\psi](x)\,\mu(\D{x})\,,
\end{align*}
with
\begin{equation*}
F[\psi](x)\,\df\, \inf_{u\in\Act}\;\sup_{w\in\Rd}\;
\bigl[\cA \psi(x,u,w) + R(x,u,w)\bigr]\,.
\end{equation*}
This underscores the discussion in the last paragraph of section~\ref{S4}.

\section{Finite state and action space}\label{S6}

For discrete time problems with finite state and action spaces
(i.e., $|S|, |U| < \infty$ in sections 2-3), one can go significantly further for
the reward maximization problem.
We recall below some results in this context from \cite{CDC}. 

 Consider a controlled Markov chain $\{Y_n\}$ on $S$ with state-dependent action
 space at state $i$ given by:
$$\tilde{U}_i \,\df\, \cup_{u \in U}(\{ u \}\times V_{i,u})\,,$$ where
$$V_{i,u} \,\df\, \biggl\{q( \cdot\,|\, i, u) \colon q( \cdot\,|\, i, u) \geq 0,
\ \sum_jq(j\,|\, i, u) = 1\biggr\}\,.$$
This is isomorphic to $\mathcal{P}(S)$.
Let $$K \,\df\, \cup_{i \in S}(\{i\}\times\tilde{U}_i)\,.$$
The (controlled) transition probabilities of $\{Y_n\}$ are
\begin{equation*}
\tilde{p}\bigl(j \,|\, i, (u,q( \cdot \,|\, i, u))\bigr) \,\df\, q(j \,|\, i,u)\,.
\end{equation*}

Define the per stage reward $\tilde{r} \colon K\times S \mapsto \mathcal{R}$ by:
\begin{equation*}
\tilde{r}\bigl(i, (u,q(\cdot \,|\, i,u)), j\bigr) \,\df\,
r(i, u, j) - D\bigl(q(\cdot \,|\, i, u) \|\, p(\cdot\,|\, i, u)\bigr)\,.
\end{equation*}
Let $\{(Z_n, Q_n), n \geq 0\}$ denote the $\tilde{U}_{Y_n}$-valued control process.
Consider the problem:
Maximize the long run average reward
$$\liminf_{N\uparrow\infty}\frac{1}{N}\sum_{n=0}^{N-1}
\Exp\left[\tilde{r}\left(Y_n, (Z_n, Q_n), Y_{n+1}\right)\right].$$

\bigskip

Define the corresponding ergodic occupation measure $\gamma \in \mathcal{P}(K\times S)$ by
\begin{equation*}
\gamma(i, (u, \D q), j) \,\df\, \gamma_1(i)\gamma_2(u,\D q\,|\, i)
\gamma_3\bigl(j\,|\, i, (u, q)\bigr)\,,
\end{equation*}
where $\gamma_1$ is an invariant probability distribution (not necessarily unique)
under the transition kernel
\begin{equation*}
\check{\gamma}(j\,|\, i) \,=\, \sum_u\int_{V_{i,u}}\gamma_2(u,\D q\,|\, i)
\gamma_3\bigl(j\,|\, i, (u, q)\bigr)\,.
\end{equation*}
Let $\mathcal{E}$ denote the set of such $\gamma$.
The above average reward control problem is equivalent to the linear program:

\medskip
\noindent \textbf{P0} \ \ \ Maximize
$$\sum_{i,j,u}\int \gamma(i, (u,\D q), j)\tilde{r}(i, (u,q), j)$$
over $\mathcal{E}$.

Recall that $\mathcal{E}$ is specified by linear constraints
and its extreme points correspond to stationary Markov policies (\cite{BorkarMC},
Chapter~V).
The maximum will be attained at an extreme point of $\mathcal{E}$ corresponding
to a stationary Markov policy.
This LP can be simplified as:

Maximize
$$\sum_{i,j}\int \gamma'(i,u,j)\Bigl[r(i, u, j)
- D\bigl(q(\cdot \,|\, i, u) \|\, p(\cdot \,|\, i, u)\bigr)\Bigr]$$
over
\begin{multline*}
\tilde{\mathcal{E}} \,\df\, \biggl\{\gamma' \in \mathcal{P}(S\times U\times S)
\colon \gamma'(i, u, j) = \gamma_1(i)\varphi(u\,|\, i)q(j \,|\, i, u),
\text{\ where\ } \gamma_1(\cdot) \text{\ is invariant}\\ \text{under the transition kernel\ }
\breve{\gamma}(j|\,i) \df \sum_u\varphi(u\,|\,i)q(j\,|\,i,u)\biggr\}\,.
\end{multline*}

\bigskip

The dual LP is:

\medskip

\noindent Minimize $\breve{\lambda}$ subject to
\begin{align*}
\breve{\lambda} &\,\ge\, \lambda(i)\,, \\
\lambda(i) + V(i) &\,\ge\, \sum_jq(j\,|\,i,u)
\bigl(\tilde{r}(i, (u, q( \cdot\,|\, i,u)), j) + V(j)\bigr)\,, \\
\lambda(i) &\,\ge\, \sum_jq(j\,|\, i,u)\lambda(j)\,, \\
&\mspace{60mu} \forall \ i \in S, \ (u, q(\cdot\,|\, i,u)) \in \tilde{U}_i\,.
\end{align*}
The proof goes through finite approximations.
Note that the LP has infinitely many constraints.
However, it does pave the way for the corresponding dynamic programming principle.
The dynamic programming formulation equivalent to the above LP turns out to be as follows:
\begin{align*}
\lambda^* &\,=\, \max_i\lambda(i)\,, \\
\lambda(i) + V(i) &\,=\, \max_{u, q(\cdot\,|\, i, u)}
\Bigl(\sum_jq(j\,|\,i, u)\bigl(V(j)
 + \tilde{r}(i, (u, q(\cdot\,|\, i, u), j))\bigr)\Bigr)\,, \qquad (\dagger) \\
\lambda(i) &\,=\, \max_{(u, q(\cdot\,|\,i, u)) \in B_i}\sum_jq(j\,|\, i, u)\lambda(j)\,,
\end{align*}
for all $i\in S$,
where $B_i$ is the Argmax in ($\dagger$). Once again, the proof goes through finite
 approximations. The maximization over $q$ in ($\dagger$) can be explicitly
 performed using the `Gibbs variational principle' from statistical mechanics.
For fixed $i,u,$ the maximum is attained at
\begin{equation*}
q^*(j\,|\, i, u) \,\df\, \frac{p(j\,|\,i,u)\E^{r(i, u, j) + V(j)}}
{\sum_kp(k\,|\,i,u)\E^{r(i, u, k) + V(k)}}\,.
\end{equation*}
Substitute back, setting
 $$\Phi(i) \,\df\, \E^{V(i)}, \quad \Lambda(i) \,\df\, \E^{\lambda(i)},\quad i \in S\,,$$
 and exponentiate both sides of ($\dagger$).
This leads to the multiplicative dynamic programming equations for infinite horizon
risk-sensitive reward in the general degenerate case:
\begin{align*}
\Lambda(i)\Phi(i) &\,=\, \max_u\sum_jp(j\,|\, i, u)\left(\E^{r(i, u, j)}\Phi(j)\right)\,,
\mspace{100mu}(\dagger\dagger) \\
\Lambda(i) &\,=\, \max_{u \in D_i}\sum_j\left(\frac{p(j\,|\,i,u)\E^{r(i, u, j)}\Phi(j)}
{\sum_kp(k\,|\,i,u)\E^{r(i, u, k)}\Phi(k)}\right)\Lambda(j)\,,
\end{align*}
for all $i\in S$, where $D_i$ is the Argmax in ($\dagger\dagger$). This is the analog of
the Howard--Kallenberg results for ergodic or `average reward'
control (\cite{Puterman}, Chapter 9).
Observe the occurrence of the `\emph{twisted kernel}', which sets it apart from the
average reward case.

\section{Future directions}

There are several directions left uncharted in this broad problem area.
Some of them are listed below.
\begin{enumerate}
\item There are some in-between cases that need to be analyzed, e.g., controlled
Markov chains with countably infinite state space.
Under the strong `Doeblin condition', the abstract Collatz-Wielandt formula has been
derived for these in \cite{Cavazos}. This needs to be extended to more general cases.

\item The counterpart of the dynamic programming equations derived for reducible
risk-sensitive reward processes can also be expected to hold for risk-sensitive cost
problems and is yet to be established.

\item Concrete computational schemes based on approximate concave maximization problems
is another direction worth pursuing.
\end{enumerate}

\section*{Acknowledgements} The work of A.A.\ was supported in part by
the National Science Foundation through grant DMS-1715210, and in part
the Army Research Office through grant W911NF-17-1-001.
The work of V.S.B.\ was supported by a J.\ C.\ Bose Fellowship from the Government of India.


\begin{bibdiv}
\begin{biblist}

\bib{Ananth}{article}{
      author={Anantharam, V.},
      author={Borkar, V.~S.},
       title={A variational formula for risk-sensitive reward},
        date={2017},
     journal={SIAM J. Control Optim.},
      volume={55},
      number={2},
       pages={961\ndash 988},
      review={\MR{3629428}},
}

\bib{AriAnup}{article}{
      author={Arapostathis, A.},
      author={Biswas, A.},
       title={A variational formula for risk-sensitive control of diffusions in
  $\mathbb{R}^d$},
        date={2018},
     journal={arXiv e-prints},
      volume={1810.01180},
      eprint={https://arxiv.org/abs/1810.01180},
}

\bib{AABK}{article}{
      author={Arapostathis, A.},
      author={Biswas, A.},
      author={Borkar, V.~S.},
      author={Suresh~Kumar, K.},
       title={A variational characterization of the risk-sensitive average
  reward for controlled diffusions in $\mathbb{R}^d$},
        date={2019},
     journal={arXiv e-prints},
      volume={1903.08346},
      eprint={https://arxiv.org/abs/1903.08346},
}

\bib{ABS-19}{article}{
      author={Arapostathis, Ari},
      author={Biswas, Anup},
      author={Saha, Subhamay},
       title={Strict monotonicity of principal eigenvalues of elliptic
  operators in {$\Bbb{R}^d$} and risk-sensitive control},
        date={2019},
     journal={J. Math. Pures Appl. (9)},
      volume={124},
       pages={169\ndash 219},
      review={\MR{3926044}},
}

\bib{ABG}{book}{
      author={Arapostathis, Ari},
      author={Borkar, Vivek~S.},
      author={Ghosh, Mrinal~K.},
       title={Ergodic control of diffusion processes},
      series={Encyclopedia of Mathematics and its Applications},
   publisher={Cambridge University Press, Cambridge},
        date={2012},
      volume={143},
      review={\MR{2884272}},
}

\bib{ABK}{article}{
      author={Arapostathis, Ari},
      author={Borkar, Vivek~S.},
      author={Kumar, K.~Suresh},
       title={Risk-sensitive control and an abstract {C}ollatz-{W}ielandt
  formula},
        date={2016},
     journal={J. Theoret. Probab.},
      volume={29},
      number={4},
       pages={1458\ndash 1484},
      review={\MR{3571250}},
}

\bib{Bielecki}{article}{
      author={Bielecki, Tomasz},
      author={Hern\'{a}ndez-Hern\'{a}ndez, Daniel},
      author={Pliska, Stanley~R.},
       title={Risk sensitive control of finite state {M}arkov chains in
  discrete time, with applications to portfolio management},
        date={1999},
     journal={Math. Methods Oper. Res.},
      volume={50},
      number={2},
       pages={167\ndash 188},
      review={\MR{1732397}},
}

\bib{BorkarMC}{book}{
      author={{Borkar}, V.~S.},
       title={Topics in controlled markov chains},
      series={Pitman Research Notes in Mathematics},
   publisher={Logman Scientific and Technical, Harlow},
        date={1991},
      volume={240},
      review={\MR{1619036}},
}

\bib{CDC}{inproceedings}{
      author={{Borkar}, V.~S.},
       title={Linear and dynamic programming approaches to degenerate
  risk-sensitive reward processes},
        date={2017},
   booktitle={2017 {IEEE} 56th {A}nnual {C}onference on {D}ecision and
  {C}ontrol ({CDC})},
       pages={3714\ndash 3718},
}

\bib{Cavazos}{article}{
      author={Cavazos-Cadena, Rolando},
       title={Characterization of the optimal risk-sensitive average cost in
  denumerable {M}arkov decision chains},
        date={2018},
     journal={Math. Oper. Res.},
      volume={43},
      number={3},
       pages={1025\ndash 1050},
      review={\MR{3846082}},
}

\bib{Pagter}{article}{
      author={de~Pagter, Ben},
       title={Irreducible compact operators},
        date={1986},
     journal={Math. Z.},
      volume={192},
      number={1},
       pages={149\ndash 153},
      review={\MR{835399}},
}

\bib{Dembo}{book}{
      author={Dembo, Amir},
      author={Zeitouni, Ofer},
       title={Large deviations techniques and applications},
     edition={Second},
      series={Applications of Mathematics (New York)},
   publisher={Springer-Verlag, New York},
        date={1998},
      volume={38},
      review={\MR{1619036}},
}

\bib{DoVa}{article}{
      author={Donsker, Monroe~D.},
      author={Varadhan, S. R.~S.},
       title={On a variational formula for the principal eigenvalue for
  operators with maximum principle},
        date={1975},
     journal={Proc. Nat. Acad. Sci. U.S.A.},
      volume={72},
       pages={780\ndash 783},
      review={\MR{0361998}},
}

\bib{Flem}{article}{
      author={Fleming, Wendell~H.},
      author={McEneaney, William~M.},
       title={Risk-sensitive control on an infinite time horizon},
        date={1995},
     journal={SIAM J. Control Optim.},
      volume={33},
      number={6},
       pages={1881\ndash 1915},
      review={\MR{1358100}},
}

\bib{Krein}{article}{
      author={Kre\u{\i}n, M.~G.},
      author={Rutman, M.~A.},
       title={Linear operators leaving invariant a cone in a {B}anach space},
        date={1950},
     journal={Amer. Math. Soc. Translation},
      volume={1950},
      number={26},
       pages={128},
      review={\MR{0038008}},
}

\bib{Meyer}{book}{
      author={Meyer, Carl},
       title={Matrix analysis and applied linear algebra},
   publisher={Society for Industrial and Applied Mathematics (SIAM),
  Philadelphia, PA},
        date={2000},
      review={\MR{1777382}},
}

\bib{Ogiwara}{article}{
      author={Ogiwara, Toshiko},
       title={Nonlinear {P}erron-{F}robenius problem on an ordered {B}anach
  space},
        date={1995},
     journal={Japan. J. Math. (N.S.)},
      volume={21},
      number={1},
       pages={42\ndash 103},
      review={\MR{1338356}},
}

\bib{Puterman}{book}{
      author={Puterman, Martin~L.},
       title={Markov decision processes: discrete stochastic dynamic
  programming},
      series={Wiley Series in Probability and Mathematical Statistics: Applied
  Probability and Statistics},
   publisher={John Wiley \& Sons, Inc., New York},
        date={1994},
        note={A Wiley-Interscience Publication},
      review={\MR{1270015}},
}

\end{biblist}
\end{bibdiv}

\end{document}